\documentclass[conference]{IEEEtran}
\IEEEoverridecommandlockouts
\usepackage{cite}
\usepackage{amsmath,amssymb,amsfonts}
\usepackage{algorithmic}
\usepackage{graphicx}
\usepackage{textcomp}
\usepackage{xcolor}
\usepackage{subcaption}
\usepackage{multirow}
\usepackage{longtable}

\def\BibTeX{{\rm B\kern-.05em{\sc i\kern-.025em b}\kern-.08em
    T\kern-.1667em\lower.7ex\hbox{E}\kern-.125emX}}
\begin{document}

\title{Analysis of biologically plausible neuron models for regression with spiking neural networks
\thanks{This work was supported by the Vannevar Bush Faculty Fellowship award (GEK) from ONR (N00014-22-1-2795). The work of GEK is supported by the U.S. Department of Energy, Advanced Scientific Computing Research program, under the Scalable, Efficient and Accelerated Causal Reasoning Operators, Graphs and Spikes for Earth and Embedded Systems (SEA-CROGS) project, DE-SC0023191.}
}

\author{\IEEEauthorblockN{Mario De Florio}
\IEEEauthorblockA{\textit{Division of Applied Mathematics} \\
\textit{Brown University}\\
Providence, RI \\
mario\_de\_florio@brown.edu}
\and
\IEEEauthorblockN{Adar Kahana}
\IEEEauthorblockA{\textit{Division of Applied Mathematics} \\
\textit{Brown University}\\
Providence, RI \\
adar\_kahana@brown.edu}
\and
\IEEEauthorblockN{George Em Karniadakis}
\IEEEauthorblockA{\textit{Division of Applied Mathematics} \\
\textit{Brown University}\\
Providence, RI \\
george\_karniadakis@brown.edu}

}

\maketitle

\begin{abstract}

This paper explores the impact of biologically plausible neuron models on the performance of Spiking Neural Networks (SNNs) for regression tasks. While SNNs are widely recognized for classification tasks, their application to Scientific Machine Learning and regression remains underexplored. We focus on the membrane component of SNNs, comparing four neuron models: Leaky Integrate-and-Fire, FitzHugh–Nagumo, Izhikevich, and Hodgkin-Huxley. We investigate their effect on SNN accuracy and efficiency for function regression tasks, by using Euler and Runge-Kutta 4th-order approximation schemes. We show how more biologically plausible neuron models improve the accuracy of SNNs while reducing the number of spikes in the system. The latter represents an energetic gain on actual neuromorphic chips since it directly reflects the amount of energy required for the computations.

\end{abstract}


\section{Introduction}

The neuromorphic research has been gaining a lot of traction recently. The main reasons are the potential efficiency, that unlocks faster training and inference when using neural networks \cite{markovic2020physics}. Another reason is the scalability. With more neuromorphic hardware \cite{schemmel2010wafer} and software \cite{titirsha2021role} made available every day, researchers are inclined to explore applying their methods developed on traditional hardware to neuromorphic chips. Finally, the biological plausibility of neuromorphic research is a backwind that drives this research \cite{yang2020neuromorphic}. The core of neuromorphic research is currently Spiking Neural Networks (SNNs). Several works \cite{han2020deep,lemaire2022analytical,misra2024use,putra2020fspinn,park2019fast} have shown the superiority of SNNs in terms of efficiency, while also maintaining sufficient accuracy compared to their Artificial Neural Network (ANN) counterparts. Most works so far that utilize SNNs for Machine Learning (ML) focus on classification tasks. This is due to the natural usage of neuromorphic hardware for classification, where the inference is usually binary - $0$ if the sample belongs to the specific class and $1$ otherwise (also referred to as one-hot classification). This aligns with the neuronal model well - $0$ for no spike and $1$ for spike.

A less explored area of neuromorphic research is the use of SNNs for Scientific Machine Learning (SciML). An example of a SciML task would be approximating a continuous function using ML. These problems are often using regression instead of classification, which is a challenge for SNNs since they do not naturally fit the paradigm like in classification. However, to some extent, regression is still possible using SNNs as presented in \cite{kahana2022function}. Following regression, we have shown that SNNs can be used for operator regression, numerical schemes for dynamical systems \cite{fatunla2014numerical}, physics-informed machine learning \cite{karniadakis2021physics}, and more \cite{yu2014brain,mostafa2017supervised,comsa2020temporal}. There is still a lot to study, revolving around neuromorphic algorithms for SciML, and this work presents a fundamental step.

\begin{figure*}[t]
    \centering
    \includegraphics[width=\textwidth]{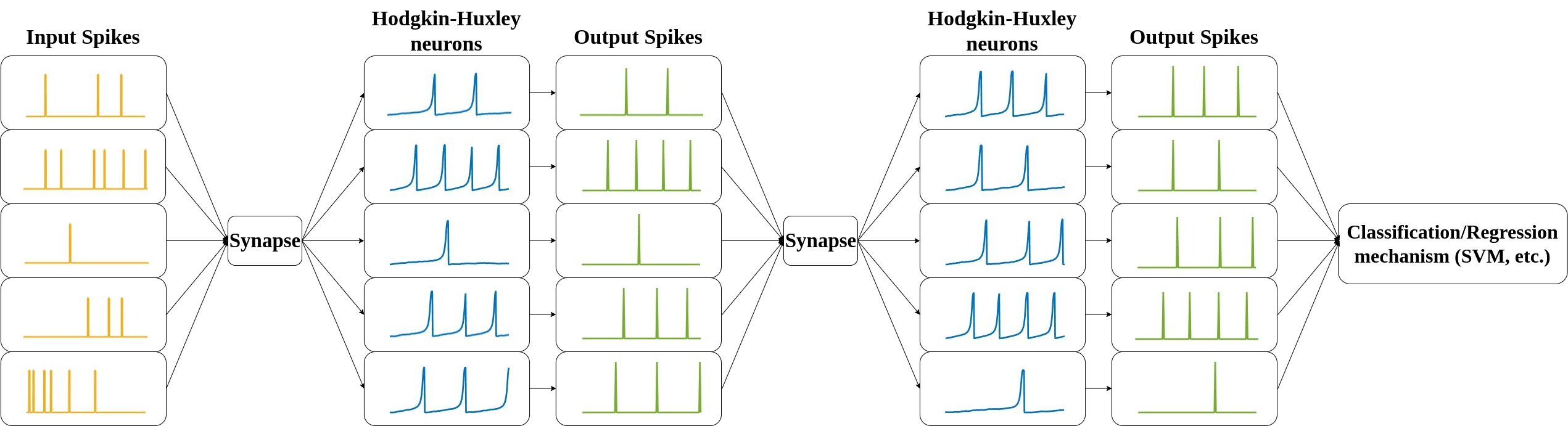}
    \caption{Sketch of an SNN architecture. The spiking data is used as inputs to Hodgkin-Huxley neurons (membranes) and then fed into a synapse. In this architecture, we have a membrane potential that creates output spikes, followed by a synapse, and another membrane that generates output spikes, followed by a second synapse. Finally, a classification or regression mechanism should be applied. In this study, we employ only one membrane followed by a single synapse to perform regression tasks, as elaborated in the results section.}
    \label{fig:SNN_scheme}
\end{figure*}

SNNs have three main components. The first is the input data, which is in the form of spikes (binary zeros and ones). Data is encoded into spikes before being input into the SNN, and outputs are often decoded back to at the end (depending on the task). The second element is the synapse, which can be thought of as a fully connected layer in the classical deep learning sense. The synapses are trained using learning rules \cite{bengio1990learning}. Last are the membranes, which are analogous to the activation functions. An illustration of the structure of an SNN is presented in Fig. \ref{fig:SNN_scheme} The membranes are the focus of this work.

From the definition of the biological neuron \cite{levitan2015neuron}, the membrane is a voltage gate at the entrance of the neuronal cell. The membrane ``decides'' how much voltage enters, and when a certain voltage threshold is crossed the neuron spikes for a short period of time, and the membrane voltage resets to its resting voltage. The first discovery of the membrane model was in 1907 by Louis Lapicque \cite{lapicque1907recherches}, where he injected a frog's leg with voltage until it twitched, inventing the Leaky-Integrate and Fire (LIF) neuron model that is commonly used in many studies. This model is famous for its simplicity, and also biological accuracy. Many neuroscience studies \cite{laing2009stochastic,ostojic2011spiking,krasovskaya2019salience,kass2018computational} use this model, and all SNN software and hardware simulators we are aware of have implementations of a LIF model. The LIF neuron is modeled using a Resistor-Capacitor circuit (hence the term neuron model), as explained in detail in the following section.

Several neuroscience-related studies have been conducted to develop more accurate membrane formulations. The FitzHugh-Nagumo model \cite{fitzhugh1961impulses,nagumo1962active} introduces a set of two coupled Ordinary Differential Equations (ODEs) to better capture the oscillatory behavior of the neuronal cell. The Izhikevich model \cite{izhikevich2003simple} is also a set of ODEs, which captures better the spiking patterns. The Hodgkin-Huxley model \cite{hodgkin1952quantitative} is a set of four non-linear coupled ODEs that goes deeper into the transfer of potassium and sodium in and out of the neuronal cell. While the aforementioned three models and LIF are explored in this work, we note that several other models have been proposed \cite{pinsky1994intrinsic,hayman1999mcculloch,lecar2007morris,storace2008hindmarsh}. The trade-off when estimating the eligibility of the models for SciML is between the accuracy of the overall process as described in the Results section and the complexity of the model, which directly impacts its efficiency. These were the guidelines for choosing the models to compare in this work, where LIF is the simplest and most efficient and Hodgkin-Huxley on the other accuracy end of the spectrum. The intermediate ones are, in addition to their qualities, very popular amongst this scientific community.\\

The main question we explore in this work is: {\em how do the different membrane models contribute to the success of an SNN?} As mentioned, SNNs are measured by their overall accuracy, as well as their efficiency. We first develop a generic implementation of a SNN with a synapse and a membrane, and simply alternate between the different membrane models. To do this, we implemented all neuron models and fine-tuned the parameters to achieve the desired results. The findings we report in this manuscript bring valuable insights for those experimenting with SNNs. In addition, we find that more complex models, although taking longer to compute, reduce the total number of spikes registered, which means lower overall energy (higher efficiency).

\section{Biological Neuron Models}

\subsection{Leaky Integrate-and-Fire model} 

The Leaky Integrate-and-Fire (LIF) model is a simplified neuronal model widely used for its computational efficiency and ease of analysis. It provides a basic representation of neuronal dynamics, capturing the essential features of membrane potential integration and action potential generation. The earliest form of the Integrate-and-Fire model was originally introduced by Lapicque in 1907 \cite{lapicque1907recherches}, and subsequently, the LIF model was developed after adding an attenuation term \cite{abbott1999lapicque}.\\
The dynamics of the LIF model is governed by a single ODE that describes the evolution of the membrane potential $V(t)$ over time:
\begin{equation}
    \dfrac{d V(t)}{dt} = \frac{V(t) - V_{rest}}{RC} + I(t),
    \label{eq:LIF_model}
\end{equation}
where $C$ is the membrane capacitance, $R$ is the membrane resistance, $I(t)$ is the input current, and $V_{rest}$ is the membrane resting voltage. In our simulation, $RC = 0.2$ and $I_c(t)$ are input currents of magnitude $1$.\\
The LIF model is commonly used to study the response of a neuron to input currents. A spike is generated when the membrane potential reaches a threshold value $V_{\text{th}}=1$. Upon reaching this threshold, the membrane potential is reset to a resting value $V_{rest}$. The thresholding mechanism is expressed as follows:
\begin{equation}
\text{if } V \geq V_{th}, \text{ then spike and } V \rightarrow V_{rest}
\label{eq:LIF_threshold}
\end{equation}

\subsection{FitzHugh–Nagumo model}

The FitzHugh–Nagumo (FHN) model is a simplified two-variable system of ODEs that was introduced by Richard FitzHugh in 1961 \cite{fitzhugh1961impulses} and later refined by J. Nagumo et al. in 1962 \cite{nagumo1962active}. The FHN model serves as a mathematical abstraction to capture the characteristics of neuronal excitability and oscillatory behavior, and it has the following form
\begin{equation}
    \begin{cases}
        \dfrac{dV}{dt} = V - \dfrac{V^3}{3} - W + I \\
        \dfrac{dW}{dt} = \dfrac{1}{\gamma} (V + \alpha - \beta W)
    \end{cases}  \;  \text{s.t.}  \qquad
    \begin{cases}
        V(0) = 0 \\
         W(0) = 0 
    \end{cases}
    \label{eq:FHN_model}
\end{equation}
where $V$ is the membrane potential, $W$ is a recovery variable, and $I$ is the input current. The parameters governing the system's behavior are the bifurcation parameter $\alpha = 0.7$ determining the onset of excitability, $\beta=0.8$ which controls the recovery variable's sensitivity to the membrane potential, and $\gamma=12.5$ governs the timescale of the recovery variable. A spike is generated and the membrane resets to the resting value 0 once the threshold value of 1 is reached.

\subsection{Izhikevich model}

The Izhikevich model, introduced by Eugene M. Izhikevich in 2003 \cite{izhikevich2003simple}, is a simple yet powerful two-variable neuronal model that captures a wide range of spiking behaviors observed in real neurons. This model is particularly valuable for its computational efficiency and ability to replicate diverse neural firing patterns.\\
The dynamics of the Izhikevich model are governed by the following set of equations:

\begin{equation}
    \begin{cases}
        \dfrac{dv}{dt} = 0.04v^2 + 5v + 140 - u + I \\
        \dfrac{du}{dt} = a(bv - u)
    \end{cases}  \;  \text{s.t.}  \quad
    \begin{cases}
        v(0) = -70 mV \\
         u(0) = -14
    \end{cases}
    \label{eq:Izhikevich_model}
\end{equation}
with the auxiliary condition
\begin{equation}
\text{if } v \geq v_{th}, \text{ then spike and } \begin{cases}
    v \rightarrow c \\ u \rightarrow u + d
\end{cases} 
\label{eq:izh_threshold}
\end{equation}
where $v$ represents the membrane potential, $u$ is a recovery variable, and $I=10mA$ is the input current. The parameters shaping the spiking dynamics are $a=0.02 \frac{1}{ms}$, which controls the time scale of the recovery variable, $b=0.2$ that influences the sensitivity of the recovery variable to the membrane potential, $c=-50 mV$, which is the after-spike reset value of the membrane potential caused by the fast high-threshold $K^+$ conductances, and $d=2mV$, which is the after-spike reset of the recovery variable caused by slow high-threshold $Na^+$ and $K^+$ conductances.

\subsection{Hodgkin-Huxley model}

\begin{figure}[t]
    \centering
\includegraphics[width=\linewidth]{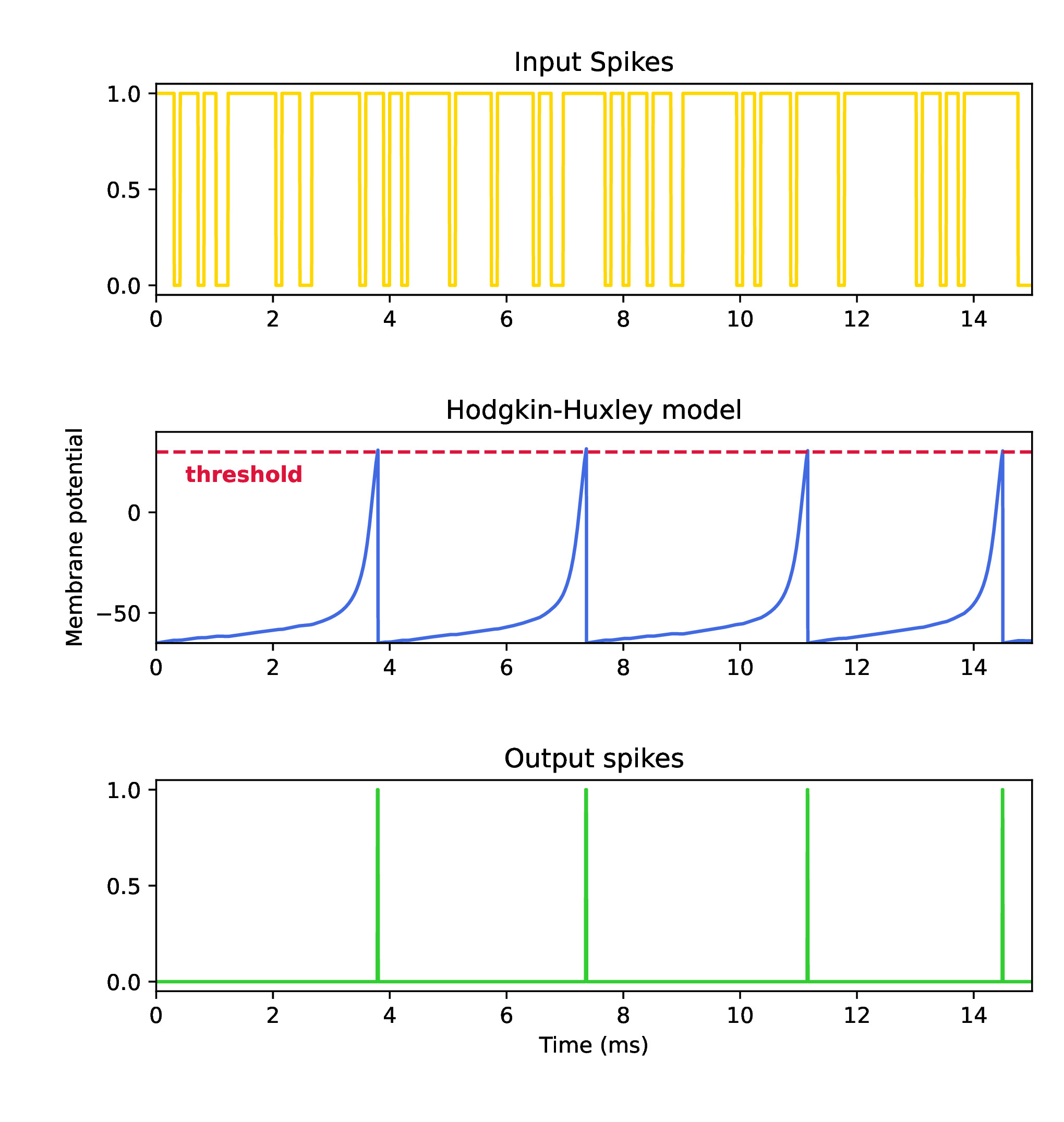}
    \caption{Example of the membrane potential over time for given input spikes. The reset mechanism is set to a resting position without refraction. The top plot presents the input spikes over time. The middle plot presents the membrane activity modeled by the Hodgkin-Huxley model. Notice that when the membrane potential reaches the threshold, an output spike is registered, and the membrane resets. The bottom plot presents the registered output spikes over time.}
    \label{fig:membrane}
\end{figure}

The Hodgkin-Huxley model, proposed in 1952 by British physiologists Alan Lloyd Hodgkin and Andrew Huxley \cite{hodgkin1952quantitative}, revolutionized our understanding of the mechanisms underlying the generation and propagation of action potentials in neurons. This mathematical model provides a framework for describing the complex dynamics of ion channels and membrane potentials in excitable cells, particularly neurons. The activation of the membrane potential is described by the following system of four non-linear Ordinary Differential Equations (ODEs) subject to certain initial conditions \cite{teka2016power}:
\begin{equation}
    \begin{cases}
      C_m \dfrac{d V_m}{dt} =  - g_l(V_m - e_l) - g_kn^4(V_m - e_k)\\
      \qquad \qquad - g_{Na} m^3 h(V_m - e_{Na}) + I \\
     \dfrac{d n}{dt}  = \alpha_n(V_m) (1 - n) - \beta_n(V_m) n\\
      \dfrac{d m}{dt}  = \alpha_m (V_m)(1 - m) - \beta_m(V_m) m\\
      \dfrac{d h}{dt}  = \alpha_h (V_m)(1 - h) - \beta_h(V_m) h\\
    \end{cases} 
    \label{eq:HH_system}
\end{equation}
subject to $ V_m(0) = -65 mV, n(0) = 0.3177, m(0) = 0.0529$, and $h(0) = 0.5960$. Each ODE describes the time evolution of the four variables of the HH model. In particular, $V_m$ is the membrane potential, $n$ and $m$ are the activation gating variables for potassium and sodium channels, respectively, and $h$ is the inactivation gating variable for the sodium channels. $C_m$ represents the membrane capacitance, $g_l$, $g_K$ and $g_{Na}$ are the leak, potassium and sodium conductances, respectively, $e_l$, $e_k$ and $e_{Na}$ are their reversal potentials, and $I$ is the input current per unit area. The $\alpha$ and $\beta$ functions represent the rate constants for the ion channels, and they can be expressed as
\begin{align*}
    & \alpha_n(V_m) = \dfrac{0.1 - 0.0 1U_{m0}}{e^{1 - 0.1 U_{m0}}-1} \qquad  &  \beta_n(V_m) = 0.125  e^{- \dfrac{U_{m0}}{80} }   \\
    & \alpha_m (V_m)= \dfrac{2.5 - 0.1 U_{m0}}{e^{2.5 - 0.1 U_{m0}}-1} \qquad  & \beta_m (V_m)= 4  e^{- \dfrac{U_{m0}}{18} }  \\
    & \alpha_h (V_m)=  0.07  e^{- \dfrac{U_{m0}}{20} }  \qquad  &   \beta_h (V_m)= \dfrac{1}{1 + e^{ 3  - \dfrac{U_{m0}}{10}   } }
\end{align*}
where $ U_{m0} = V_m - V_m(0) $, with $V_m(0) = -65 mV$, and the values of the parameters used in our simulations, assuming $1cm^2$ of membrane, are $  C_m = 1 \mu F$,  $ g_l = 0.3 mS $, $ g_k = 36 mS $, $ g_{Na} = 120 mS $,  $ e_k = -77 mV $, $ e_l = -54 mV  $, $ e_{Na} = 50 mV $. In our simulations, the input current $I$ is represented by spikes of $5mA$ of intensity. The membrane has a threshold parameter set at $30 mV$. When the membrane potential reaches the threshold, a spike is registered in the neuron, and the membrane resets to the value of $V_m(0)$. An example of this membrane's behavior is represented in Figure \ref{fig:membrane}.

\section{Numerical methods}

\subsection{Euler method}

To simplify the calculations, we can employ the forward Euler finite difference method to directly approximate the solution of the Hodgkin-Huxley model \cite{hanslien2005maximum} from equations \eqref{eq:HH_system}, as follows:
\begin{equation}
    \begin{aligned}
        V_m(t+\Delta t) &= V_m(t) + \Delta t  \bigl(I(t) - g_K n^4  (V_m (t) - e_K) \\ 
        - g_{Na} m^3 &  h  (V_m(t) - e_{Na}) - g_L  (V_m(t) - e_L)\bigr) / C_m  \\
        n(t+\Delta t) &= n(t) + \Delta t  (\alpha_n  (1 - n(t)) - \beta_n  n(t)) \\
        m(t+\Delta t) &= m(t) + \Delta t  (\alpha_m  (1 - m(t)) - \beta_m  m(t)) \\
        h(t+\Delta t) &= h(t) + \Delta t (\alpha_h  (1 - h(t)) - \beta_h  h(t)) \\
    \end{aligned}
\end{equation}
where $\Delta t$ is set equal to 0.1. This serves as a simplified approximation to the solution of the Hodgkin-Huxley model, and it's anticipated that the accuracy in computing the membrane potential will be reduced. Nevertheless, this approach offers a more lightweight implementation, enabling more efficient computations.

\subsection{Runge-Kutta 4th order method}

In contrast to the Euler method, the Runge-Kutta 4th Order (RK4) numerical integration scheme provides a more accurate approximation of the solution to the Hodgkin-Huxley model \cite{moore1974numerical}. The RK4 method involves multiple steps to iteratively update the variables, providing enhanced accuracy in capturing the intricate dynamics of the membrane potential and gating variables.\\
The slopes $k_1, k_2, k_3$, and $k_4$ for all the HH variables are
\begin{equation}
\begin{aligned}
k_{1V} &= \Delta t \bigl(I - g_K n^4(V_m - e_K) - g_{Na} m^3h(V_m - e_{Na}) \\
&- g_L(V_m - e_L)\bigr) / C_m  \\
k_{1n} &= \Delta t (\alpha_n(1 - n) - \beta_n n) \\
k_{1m} &= \Delta t (\alpha_m(1 - m) - \beta_m m) \\
k_{1h} &= \Delta t (\alpha_h(1 - h) - \beta_h h)
\end{aligned}
\end{equation}

\begin{equation}
\begin{aligned}
&k_{2V} = \Delta t \bigl(I - g_K (n + 0.5k_{1n})^4(V_m + 0.5k_{1V} - e_K)\\ 
&  - g_{Na} (m + 0.5k_{1m})^3(h + 0.5k_{1h})(V_m + 0.5k_{1V} - e_{Na}) \\ 
& \qquad  - g_L (V_m + 0.5k_{1V} - e_L)\bigr) / C_m \\
&k_{2n} = \Delta t (\alpha_n(1 - (n + 0.5k_{1n})) - \beta_n (n + 0.5k_{1n})) \\
&k_{2m} = \Delta t (\alpha_m(1 - (m + 0.5k_{1m})) - \beta_m (m + 0.5k_{1m})) \\
&k_{2h} = \Delta t (\alpha_h(1 - (h + 0.5k_{1h})) - \beta_h (h + 0.5k_{1h}))
\end{aligned}
\end{equation}

\begin{equation}
\begin{aligned}
&k_{3V} = \Delta t \bigl(I - g_K (n + 0.5k_{2n})^4(V_m+ 0.5k_{2V} - e_K) \\
&  - g_{Na} (m + 0.5k_{2m})^3(h + 0.5k_{2h})(V_m+ 0.5k_{2V} - e_{Na})\\
& \qquad  - g_L (V_m+ 0.5k_{2V} - e_L)\bigr) / C_m \\
&k_{3n} = \Delta t (\alpha_n(1 - (n + 0.5k_{2n})) - \beta_n (n + 0.5k_{2n})) \\
&k_{3m} = \Delta t (\alpha_m(1 - (m + 0.5k_{2m})) - \beta_m (m + 0.5k_{2m})) \\
&k_{3h} = \Delta t (\alpha_h(1 - (h + 0.5k_{2h})) - \beta_h (h + 0.5k_{2h})) 
\end{aligned}
\end{equation}

\begin{equation}
\begin{aligned}
&k_{4V} = \Delta t \bigl(I - g_K (n + k_{3n})^4(V_m+ k_{3V} - e_K)\\
& \qquad - g_{Na} (m + k_{3m})^3(h + k_{3h})(V_m+ k_{3V} - e_{Na})  \\
& \qquad  - g_L (V_m+ k_{3V} - e_L)\bigr) / C_m \\
&k_{4n} = \Delta t (\alpha_n(1 - (n + k_{3n})) - \beta_n (n + k_{3n})) \\
&k_{4m} = \Delta t (\alpha_m(1 - (m + k_{3m})) - \beta_m (m + k_{3m})) \\
&k_{4h} = \Delta t (\alpha_h(1 - (h + k_{3h})) - \beta_h (h + k_{3h}))
\end{aligned}
\end{equation}

Thus, the updated equations for the RK4 method applied to the Hodgkin-Huxley model are as follows:
\begin{equation}
\begin{aligned}
&V_m(t + \Delta t) = V_m+ \frac{1}{6}(k_{1V} + 2k_{2V} + 2k_{3V} + k_{4V}) \\
&n(t + \Delta t) = n + \frac{1}{6}(k_{1n} + 2k_{2n} + 2k_{3n} + k_{4n}) \\
&m(t + \Delta t) = m + \frac{1}{6}(k_{1m} + 2k_{2m} + 2k_{3m} + k_{4m}) \\
&h(t + \Delta t) = h + \frac{1}{6}(k_{1h} + 2k_{2h} + 2k_{3h} + k_{4h}) 
\end{aligned}
\end{equation}
Here, $\Delta t$ is the time step, and $V_m, n, m,$ and $h$ represent the membrane potential and gating variables. The RK4 method ensures a more accurate numerical approximation of the Hodgkin-Huxley model's dynamics compared to the Euler method, making it suitable for simulations that require higher precision.

\section{Results and Discussions}

To test the performance of the different neural models computed with the Euler and RK4 methods we follow the experimental setup presented in \cite{kahana2022function} for function regression. We use the same codes, changing only the activation from LIF neurons into the various ones presented in this work. Since we do not have access to a neuromorphic chip for testing the method with spikes, we are using a GPU-based training with TensorFlow 2 \cite{raschka2019python}. The energetic gains are measured by the number of spikes in the system, rather than execution time. This is a valid measure, as on an actual chip the number of spikes directly reflects the amount of energy required for the computations. 

\begin{figure*}[!htb]
    \centering
    \begin{subfigure}{0.32\textwidth}
        \centering
        \includegraphics[width=\linewidth]{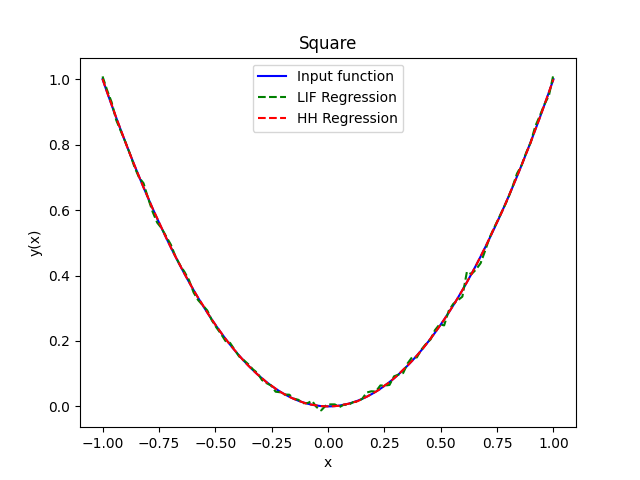}
    \end{subfigure}
    \hfill
    \begin{subfigure}{0.32\textwidth}
        \centering
        \includegraphics[width=\linewidth]{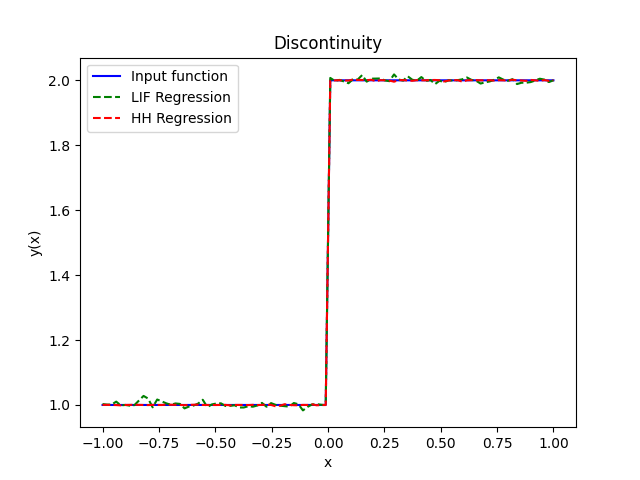}
    \end{subfigure}
    \hfill
    \begin{subfigure}{0.32\textwidth}
        \centering
        \includegraphics[width=\linewidth]{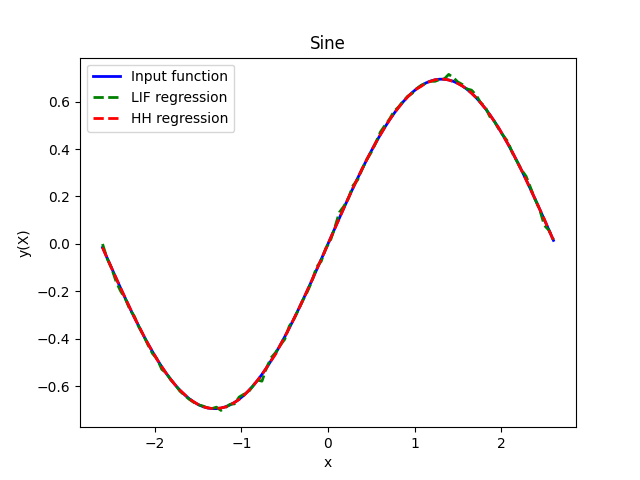}
    \end{subfigure}

    \begin{subfigure}{0.32\textwidth}
        \centering
        \includegraphics[width=\linewidth]{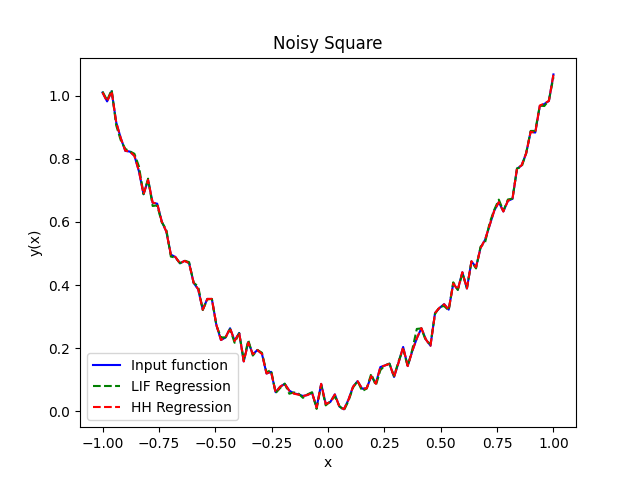}
    \end{subfigure}
    \hfill
    \begin{subfigure}{0.32\textwidth}
        \centering
        \includegraphics[width=\linewidth]{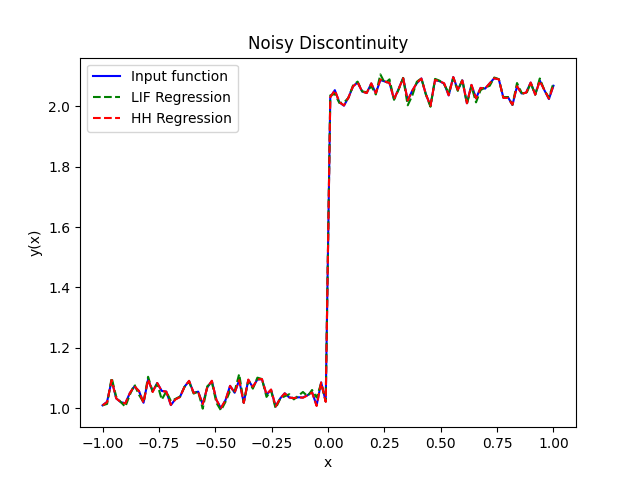}
    \end{subfigure}
    \hfill
    \begin{subfigure}{0.32\textwidth}
        \centering
        \includegraphics[width=\linewidth]{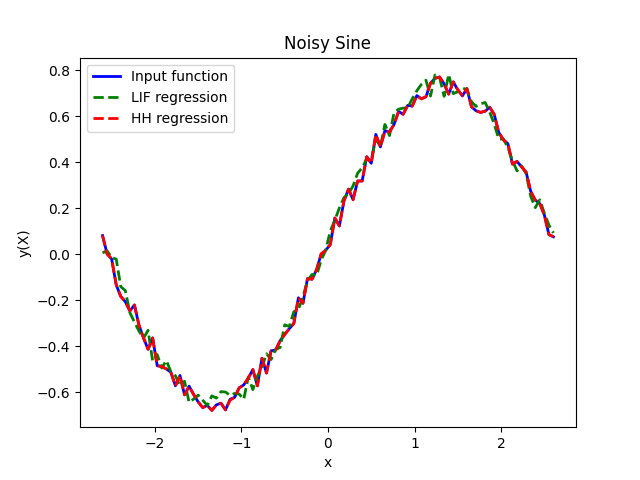}
    \end{subfigure}
    
    \caption{Comparison of function regressions using LIF model and Hodgkin-Huxley model with Euler method in the SNN implementation for three different functions, with and without noise.}
    \label{fig:plots}
\end{figure*}

\begin{table*}[!htb]
\centering
\caption{Relative $L_2$ error and number of output spikes for neuron membrane models and solvers, for both noiseless and noisy functions. The encoded neuron models are Leaky Integrate-and-Fire (LIF), FitzHugh–Nagumo (FHN), Izhikevich (IZH), and Hodgkin-Huxley (HH).}\label{tab:results}
\begin{tabular*}{\textwidth}{@{\extracolsep{\fill}}ccccccccc}
\hline \hline
   &       & \multicolumn{7}{c}{\textbf{Relative $L_2$ error (\# of output spikes)}}    \\ \hline
\multicolumn{1}{c|}{\multirow{2}{*}{\textbf{\begin{tabular}[c]{@{}c@{}}Neuron\\ model\end{tabular}}}} & \multicolumn{1}{c|}{\multirow{2}{*}{\textbf{\begin{tabular}[c]{@{}c@{}}Numerical\\ method\end{tabular}}}} & \multicolumn{3}{c|}{\textbf{Noiseless}}     & \multicolumn{3}{c|}{\textbf{Noisy}}            & \multirow{2}{*}{\textbf{\begin{tabular}[c]{@{}c@{}}Time\\ (s)\end{tabular}}} \\ \cline{3-8}
\multicolumn{1}{c|}{}         & \multicolumn{1}{c|}{}       & \textbf{Discontinuity} & \textbf{Square} & \multicolumn{1}{c|}{\textbf{Sine}} & \textbf{Discontinuity} & \textbf{Square} & \multicolumn{1}{c|}{\textbf{Sine}} &     \\ \hline \hline
\multicolumn{1}{c|}{\multirow{4}{*}{LIF}}                                                             & \multicolumn{1}{c|}{\multirow{2}{*}{Euler}}                                                               & 4.50e-03       & 1.30e-02        & \multicolumn{1}{c|}{1.43e-02}      & 4.30e-03       & 1.18e-02        & \multicolumn{1}{c|}{1.52e-02}      & \multirow{2}{*}{12}                                                          \\
\multicolumn{1}{c|}{}                                                                                 & \multicolumn{1}{c|}{}                                                                                     & (2271)         & (2242)          & \multicolumn{1}{c|}{(2312)}        & (2257)         & (2212)          & \multicolumn{1}{c|}{(2225)}        &                                                                              \\ \cline{2-9} 
\multicolumn{1}{c|}{}                                                                                 & \multicolumn{1}{c|}{\multirow{2}{*}{RK}}                                                                  & 1.864e-03      & 6.13e-03        & \multicolumn{1}{c|}{5.82-03}       & 1.66e-03       & 4.78e-03        & \multicolumn{1}{c|}{5.53e-03}      & \multirow{2}{*}{22}                                                          \\
\multicolumn{1}{c|}{}                                                                                 & \multicolumn{1}{c|}{}                                                                                     & (904)          & (886)           & \multicolumn{1}{c|}{(896)}         & (908)          & (902)           & \multicolumn{1}{c|}{(897)}         &                                                                              \\ \hline
\multicolumn{1}{c|}{\multirow{4}{*}{FHN}}                                                             & \multicolumn{1}{c|}{\multirow{2}{*}{Euler}}                                                               & 6.63e-04       & 1.96e-03        & \multicolumn{1}{c|}{1.87e-03}      & 5.84e-04       & 2.03e-03        & \multicolumn{1}{c|}{2.19e-03}      & \multirow{2}{*}{16}                                                          \\
\multicolumn{1}{c|}{}                                                                                 & \multicolumn{1}{c|}{}                                                                                     & (281)          & (269)           & \multicolumn{1}{c|}{(271)}         & (266)          & (276)           & \multicolumn{1}{c|}{(272)}         &                                                                              \\ \cline{2-9} 
\multicolumn{1}{c|}{}                                                                                 & \multicolumn{1}{c|}{\multirow{2}{*}{RK}}                                                                  & 6.16e-04       & 5.49e-04        & \multicolumn{1}{c|}{6.10e-04}      & 6.24e-04       & 5.02e-04        & \multicolumn{1}{c|}{6.40e-04}      & \multirow{2}{*}{48}                                                          \\
\multicolumn{1}{c|}{}                                                                                 & \multicolumn{1}{c|}{}                                                                                     & (55)           & (53)            & \multicolumn{1}{c|}{(51)}          & (58)           & (54)            & \multicolumn{1}{c|}{(56)}          &                                                                              \\ \hline
\multicolumn{1}{c|}{\multirow{4}{*}{IZH}}                                                             & \multicolumn{1}{c|}{\multirow{2}{*}{Euler}}                                                               & 5.78e-04       & 2.20e-03        & \multicolumn{1}{c|}{2.05e-03}      & 6.58e-04       & 1.97e-03        & \multicolumn{1}{c|}{1.91e-03}      & \multirow{2}{*}{20}                                                          \\
\multicolumn{1}{c|}{}                                                                                 & \multicolumn{1}{c|}{}                                                                                     & (219)          & (222)           & \multicolumn{1}{c|}{(217)}         & (216)          & (218)           & \multicolumn{1}{c|}{(221)}         &                                                                              \\ \cline{2-9} 
\multicolumn{1}{c|}{}                                                                                 & \multicolumn{1}{c|}{\multirow{2}{*}{RK}}                                                                  & 5.16e-04       & 5.67e-04        & \multicolumn{1}{c|}{5.89e-04}      & 5.52e-04       & 6.11e-04        & \multicolumn{1}{c|}{7.32e-04}      & \multirow{2}{*}{55}                                                          \\
\multicolumn{1}{c|}{}                                                                                 & \multicolumn{1}{c|}{}                                                                                     & (58)           & (62)            & \multicolumn{1}{c|}{(64)}          & (65)           & (67)            & \multicolumn{1}{c|}{(66)}          &                                                                              \\ \hline
\multicolumn{1}{c|}{\multirow{4}{*}{HH}}                                                              & \multicolumn{1}{c|}{\multirow{2}{*}{Euler}}                                                               & 4.17e-04       &  1.26e-03     & \multicolumn{1}{c|}{ 1.24e-03}   & 4.61e-04    & 1.20e-03        & \multicolumn{1}{c|}{1.37e-03}      & \multirow{2}{*}{48}                                                          \\
\multicolumn{1}{c|}{}                                                                                 & \multicolumn{1}{c|}{}                                                                                     & (118)          & (114)           & \multicolumn{1}{c|}{(122)}         & (113)          & (109)           & \multicolumn{1}{c|}{(121)}         &                                                                              \\ \cline{2-9} 
\multicolumn{1}{c|}{}                                                                                 & \multicolumn{1}{c|}{\multirow{2}{*}{RK}}                                                                  & 1.69e-03              & 1.02e-03               & \multicolumn{1}{c|}{8.95e-04}             & 9.20e-04              & 8.07e-04               & \multicolumn{1}{c|}{9.04e-04}             & \multirow{2}{*}{130}                                                           \\
\multicolumn{1}{c|}{}                                                                                 & \multicolumn{1}{c|}{}                                                                                     & (71)              & (71)               & \multicolumn{1}{c|}{(62)}             & (74)              & (74)               & \multicolumn{1}{c|}{(86)}             &                                                                              \\ \hline \hline
\end{tabular*}

\end{table*}

The SNN architecture we use is a membrane layer, followed by a synapse. Since we are not using a neuromorphic chip, we use the back-propagation algorithm to train the synapse and learn the weights of the network. The shallow network we choose here is sufficient for demonstrating the advantages of using different membrane models and analyzing the results. We point out that using a more complex network architecture will likely yield better performance. In addition, with this specific architecture, the back-propagation algorithm operates only on the last layer, which is considered more ``biologically plausible''. The input data are encoded using $N_t = 150$ time steps and $N_x = 100$ collocation points for spatial discretization. The outputs are floating numbers (no decoding required). The three functions we attempt for regression, with and without noise, are:
\begin{itemize}
   \item $Square$
   \begin{equation}
       y(x) = x^2
   \end{equation}
    \item $Discontinuity$
    \begin{equation}
        y(x) = \begin{cases}
            1 \qquad \text{for} \qquad x \leq 0\\
            2 \qquad \text{for} \qquad x > 0
        \end{cases}
    \end{equation}
    \item $Sine$
    \begin{equation}
        y(x) = \frac{\sin(kx)}{k^2} \qquad \text{with} \qquad k =1.2
    \end{equation}
\end{itemize}
The square function is more of a sanity check. The discontinuous function is challenging, since it has a sharp jump (not differentiable), and the sine oscillates around 0 which often causes issues. We also consider the cases in which the functions are perturbed with 0.1 additive normally distributed noise. The performance of the regression is reported in terms of the $L_2$-norm difference between the predicted and exact solution as
\begin{equation}
    L_2\text{-norm} = \sum_{j=1}^{N_x} \left(y_{pred}(x_j) - y_{exact}(x_j)\right)^2,
\end{equation}
where $N_x$ is the number of collocation points of the spatial discretization. In Figure \ref{fig:plots}, a qualitative visualization of the regressions made by using the Euler method for the LIF and HH models is reported for the three functions, with and without noise. One can see how with HH there is a better overlap with the exact solutions. To quantitatively evaluate the different performances in terms of  the $L_2$-norm, efficiency (number of spikes produced), and computational times we can refer to Table \ref{tab:results}. As expected, the computational time increases from Euler to RK4, and from 1-ODE model (LIF) to 4-ODE model (HH), with a decrease in the $L_2$-norm error. One can see that in terms of accuracy, IZH and HH outperform LIF and FHN neuron models for all the test cases, but with a significant cost in computational time. It is interesting to notice how a greater accuracy corresponds to a lower number of output spikes, which results in a more efficient cost for the neuromorphic hardware in regression tasks.

\end{document}